\documentclass[11pt]{amsart}
\usepackage{amssymb,latexsym,comment,url,amsmath}

\usepackage[T1]{fontenc}
\usepackage{rotating,graphicx}

\usepackage{diagbox} 
\usepackage{xcolor}

\usepackage{currvita}
\usepackage{hyperref}
\usepackage{diagbox}

\usepackage{longtable}

\usepackage{longtable}

\newcommand{\PGL}{{\operatorname{PGL}}}

\newcommand{\Fp}{{\mathbb{F}_p}}

\newcommand{\OK}{{\mathcal{O}_K}}

\newcommand{\p}{\operatorname{\mathfrak p}}

\newcommand{\Q}{{\mathbb Q}}

\newcommand{\Orb}{\operatorname{Orb}}

\newcommand{\Z}{{\mathbb Z}}

\newcommand{\lcm}{\operatorname{lcm}}

\newenvironment{Proof}{\par\noindent{\sc Proof:}}%
                      {\hspace*{\fill}\nobreak$\Box$\par\medskip}
                       {\hspace*{\fill}\nobreak$\Box$\par\medskip}

\newtheorem{Proposition}{Proposition}[section]
\newtheorem{Theorem}[Proposition]{Theorem}
\newtheorem{Lemma}[Proposition]{Lemma}
\newtheorem{Corollary}[Proposition]{Corollary}
\newtheorem{Example}[Proposition]{Example}

\theoremstyle{definition}
\newtheorem{Definition}[Proposition]{Definition}
\newtheorem{Remark}[Proposition]{Remark}

\renewcommand{\baselinestretch}{1.1}

\begin{document}

	\title[Arithmetic progressions in polynomial orbits]{Arithmetic progressions in polynomial orbits}

	\author[M. Sadek]%
	{Mohammad~Sadek}
	\address{Faculty of Engineering and Natural Sciences, Sabanc{\i} University, Tuzla, \.{I}stanbul, 34956 Turkey}
	\email{mohammad.sadek@sabanciuniv.edu}

 \author[M. Wafik]%
 {Mohamed Wafik}
 \address{Department of Mathematics, 1523 Greene Street, LeConte College, University of South Carolina, Columbia SC, 29208, USA}
 \email{mwafik@sc.edu}

	\author[T. Yes\.{I}n]%
	{Tu\u{g}ba Yes\.{I}n}
\address{Faculty of Engineering and Natural Sciences, Sabanc{\i} University, Tuzla, \.{I}stanbul, 34956 Turkey}
 \email{tugbayesin@sabanciuniv.edu}

	\date{}

	\begin{abstract}
 Let $f$ be a polynomial with integer coefficients whose degree is at least 2. We consider the problem of covering the orbit $\Orb_f(t)=\{t,f(t),f(f(t)),\cdots\}$, where $t$ is an integer, using arithmetic progressions each of which contains $t$. Fixing an integer $k\ge 2$, we prove that it is impossible to cover $\Orb_f(t)$ using $k$ such arithmetic progressions unless $\Orb_f(t)$ is contained in one of these  progressions. In fact, we show that the relative density of terms covered by $k$ such arithmetic progressions in $\Orb_f(t)$ is uniformly bounded from above by a bound that depends solely on $k$. In addition, the latter relative density can be made as close as desired to $1$ by an appropriate choice of $k$ arithmetic progressions containing $t$ if $k$ is allowed to be large enough. 
	\end{abstract}
	
	\maketitle
\let\thefootnote\relax

	\section{Introduction}
 
  A dynamical system is a self-map $ f: S\longrightarrow S $ on a set $ S $ that allows iteration. The $ m$-${\textrm{th}} $ iterate of $ f $ is defined recursively by $ f^{0}(x)=x $ and $ f^{m}(x)=f(f^{m-1}(x)) $ when $ m\geq 1 $. The \textit{orbit} of a point $ P\in S $ under $ f $ is given by 
	\begin{equation*}
		\Orb_{f}(P)=\{f^{i}(P) : i=0,1,2\dots \}.
	\end{equation*}
If $\Orb_f(P)$ is infinite, $P$ is called a wandering point for $f$; otherwise, $P$ is called a preperiodic point for $f$. A preperiodic point $P \in S$ is said to be {\em periodic} if there exists an integer $n > 0$ such that $f^n (P)= P$, where $n$ is called the period of $P$. If $n$ is the smallest such integer, we say that $n$ is the exact period of $P$ under $f$.

The theme of studying intersections of polynomial orbits
is motivated by the Dynamical Mordell-Lang Conjecture, \cite[Conjecture 1.5.0.1]{Bell}, which states that if $\Phi$ is an endomorphism of a
quasiprojective variety $X$ over $\mathbb{C}$, $\alpha$ is any point in $X(\mathbb C)$ and $V\subset X$ is any subvariety, then
$\{n \ge  0 : \Phi^n(\alpha) \in V (\mathbb C)\}$ is a union of finitely many arithmetic progressions. The interested reader may consult \cite{Bell} for further discussions around this conjecture.
In \cite{Zieve1,Zieve2}, it was proved that if two nonlinear
complex polynomials have orbits with infinite intersection, then the polynomials must have a common iterate.

In this article, we concern ourselves with the arithmetic set-up of intersection of orbits of polynomials, see \cite{Burcu,Hindes,SadekYesin} for such arithmetic questions on the intersection of polynomial orbits. We show that the only situation when two rational polynomials $f$ and $g$, $\max(\deg f,\deg g)\ge 2$, have no common iterate and, yet, possess orbits with infinite intersection is only when exactly one of the polynomials $f$ and $g$ is linear and monic. The latter phenomenon leads us to studying the frequency of occurence of arithmetic progressions in nonlinear polynomial orbits. 

Given an integer $n$, it is easily seen that one can find finitely many congruences such that $n$ must satisfy at least one of them. A popular example is the one given by the following congruences
 $$n \equiv 0 \mod 2,\quad n \equiv 0 \mod 3,\quad n \equiv 1 \mod 4, \quad n \equiv 1 \mod 6, \quad n \equiv 11 \mod 12.$$

A system of residue classes is called a {\em cover}
of $S\subseteq \Z$ if any integer in $S$ belongs to one of these residue classes. The concept of 
covers of $\Z$ was introduced by Erd\"{o}s who conjectured that
 for every $N$ there is a covering system with distinct moduli greater than $N$, see \cite{Erdos}. The conjecture was disproved in \cite{Ho}.
 In this work, given a polynomial $f$ with integer coefficients whose degree is at least $2$ and an integer wandering point $t$, we study covers of $\Orb_f(t)$.

 We approach the question as follows. We first find necessary and sufficient conditions under which a certain iterate $f^n(t)$ satisfies a certain congruence. It turns out that the latter occurrence is closely related to the existence of primitive divisors of the members of the sequence $f^n(t)-t$, $n\ge 1$. See \cite{DoerksenHaensch+2012+465+472,abc, Rafe_jones,Holly_Krieger,Rice,SadekWafik} for discussions on primitive prime divisors and their densities in certain polynomial orbits.
 In fact, we give a complete classification of linear polynomials whose either forward or backward orbits intersect $\Orb_f(t)$ infinitely often. In addition, we compute the relative density of the intersection in $\Orb_f(t)$.

Given a system of residue classes $A=\{A_1,\cdots,A_k\}$ such that $t$ satisfies each congruence $A_i$, we show that $A$ covers $\Orb_f(t)$ if and only if there is $i_0$, $1\le i_0\le k$, such that $\{A_{i_0}\}$ covers $\Orb_f(t)$. In particular, we show that if none of the residue classes in $A$ covers $\Orb_f(t)$, then $\Orb_f(t)$ can never be covered by $A$, moreover, the relative density of integers in $\Orb_f(t)$ not covered by $A$ is positive.

Fixing an integer $k\ge1$, our analysis of the intersection between integer linear polynomial orbits and orbits of higher degree integer polynomials $f$ provides explicit families of positive rationals $\delta_f<1$ that appear as the relative density of elements in $\Orb_f(t)$ that can be covered by exactly $k$ residue classes. Furthermore, we show that $\delta_f$ is uniformly bounded from above by a constant $C_k<1$ that depends only on $k$. In other words, given any integer polynomial $f(x)$ with degree at least $2$, the relative density of integers that can be covered using $k$ congruences in $\Orb_f(t)$ cannot exceed $C_k$. Finally, we prove that as much as desired of $\Orb_f(t)$ can be covered using $k$ residue classes if $k$ is large enough.

\subsection*{Acknowledgments} The authors thank the anonymous referee for comments and suggestions that improved the manuscript. This project is supported by The Scientific and Technological Research Council of Turkey, T\"{U}B\.{I}TAK; research grant: ARDEB 1001/120F308. M. Sadek is partially funded by BAGEP Award of the Science Academy, Turkey.

\section{Intersection of polynomial orbits with linear polynomial orbits}

Throughout this work $K$ is a number field with algebraic closure $\overline K$ and ring of integers $\mathcal{O}_K$.

We recall that the $n$-th iteration of a polynomial $f(x)$ is defined to be $f^n(x)=f(f^{n-1}(x)),\, n\ge 1$, and $f^0(x)=x.$ Given $a\in \overline K$, the orbit, or forward orbit, of $a$ under $f$ is the set of images $\Orb_f(a)=\{f^n(a),n\ge 0\}.$ The backward orbit of $a$ under $f$ is the collection of its inverse images under the iterates of $f$. We can also denote by $\Orb_f^\pm(a):=\{f^n(a),n\in\Z\}$ the union of both the forward and backward orbits of a point $a$ under the iterates of $f$. A point $a\in K$ is called \textit{preperiodic} under $f$ of type $(m,n)$ if $f^{m+n}(a)=f^m(a)$ for some $m\geq 0, n\geq 1$. A point $a \in K$ is called {\em periodic} under $f$ if $a$ is preperiodic of type $(0,n)$. Moreover, if $n$ is the smallest such integer, then $a$ is said to be a periodic point of {\em exact period} $n$. 
If $a\in K$ is not preperiodic under $f$, then $a$ is called a {\em wandering} point for $f$. 

We define an equivalence relation on polynomials in $K[x]$ of a given degree $d\ge 2$ as follows. Two polynomial maps $f_1$ and $f_2$ in $K(x)$ of degree $d\ge2$ are {\em conjugate} if there is $\phi\in \PGL_2(\overline{K})$ such that $f_2=f_1^{\phi}:=\phi\circ f_1\circ \phi^{-1}$. If $\phi\in \PGL_2(K)$, then $f_1$ and $f_2$ are said to be $K$-conjugate. We remark that if $a$ is a periodic point of exact period $n$ for $f$, then $\phi(a)$ is a point of exact period $n$ for $f^{\phi}$. One can argue similarly for preperiodic points of $f$ and $f^{\phi}$. Moreover, if $f,\phi,$ and $a$ are defined over $K$ such that $f^n(a)=a$, then $g:=f^{\phi}$ and $b:=\phi(a)$ are defined over $K$ with $g^n(b)=b$.

If two complex polynomials $f$ and $g$ of degree at least $2$ have orbits with infinite intersection, then $f$
and $g$ must have a common iterate, \cite{Zieve1,Zieve2}. The assumption that both polynomials must be nonlinear is essential as may be emphasised by the example $\Orb_{2X^2+2}(1)\subset\Orb_{X+2}(0)$. 

It is obvious that if $\phi=\alpha x+\beta \in K[x]$, $\alpha\ne 0$, then $$\Orb_{f^{\phi}}(\phi(x))=\left\{ \phi(x),\phi(f(x)), \cdots, \phi(f^n(x)),\cdots\right\}$$ for any polynomial map $f\in K[x]$ and any $x\in K$. In addition, If $f(x)=ax+b$, $a\ne 0$, then $f^{\phi}(x)=ax+\alpha b-\beta(a-1)$. 

In what follows we consider the case when $g(x)=ax+b\in\Q[x]$ whereas $f(x)$ is an arbitrary polynomial. 

We remark that $$\Orb_g(x)=\{x,ax+b,a^2x+b(a+1),a^3x+b(a^2+a+1),\cdots, a^nx+b(a^n-1)/(a-1),\cdots\}.$$

We also notice that for a power of a linear map $f(x)=\frac{\beta^n}{t^{m-1}}x^m$, $n,m\in\Z$, and $g(x)=\beta x$, one has that $\Orb_f(t)\cap\Orb_g(t)$ is infinite for any $t\in \Z$. This is because $f^k(x)=\left(\frac{\beta^n}{t^{m-1}}\right)^{\frac{m^{k}-1}{m-1}}x^{m^k}$ which implies that $f^k(t)=\beta^{\frac{n(m^k-1)}{m-1}}t=g^{\frac{n(m^k-1)}{m-1}}(t)$.

We now consider the latter intersection $\Orb_f(t)\cap\Orb_g(t)$ when $f(x)$ is not a conjugate of a power of a linear polynomial.

\begin{Proposition}\label{prop:linear}
Let $f(x)\in K[x]$ be of degree at least $2$. 
Let $g(x)=ax+b\in K[x]$ be such that $\Orb_f(s)\cap \Orb_g^\pm(t)$ is infinite for some fixed $s,t\in K$. Then either $a$ is a root of unity in $\mathcal{O}_K$ or $f^\phi$ is a power of a linear polynomial, where $\phi(x) = x+b/(a-1)$.
\end{Proposition}
\begin{Proof}
Using a conjugation via $\phi(x)=x+b/(a-1)$ on both orbits, one may assume without loss of generality that $g(x)=ax$ and $t\ne 0$. We write $$f(x)=a_dx^d+\cdots+a_0=a_d\prod_{i=1}^e(x-c_i)^{r_i}, r_i\ge 1, \,a_d\ne 0,\, c_i\ne c_j\textrm{ for }i\ne j.$$ We also set $$ S=\{\p\in\OK \textrm{ is prime }: \nu_{\p}(s)<0\textrm{ or }\nu_{\p}(t)<0 \textrm{ or }\nu_{\p}(a_i)<0\textrm{ for some }0\le i\le d \},$$ where $\nu_{\p}$ is the associated discrete valuation to the prime $\p$.  

Now, we assume that $a$ is not a root of unity. The hypothesis implies that $f^m(s)= a^n t $ for infinitely many pairs of integers $(m,n)$. We also notice that for two such pairs $(m_1,n_1)$ and $(m_2,n_2)$, both $m_1\ne m_2$ and $n_1\ne n_2$, since otherwise, $f^{m_1}(s)=f^{m_2}(s)$ or $g^{n_1}(t)=g^{n_2}(t)$ respectively. The latter implies that $s$ is a preperiodic point for $f$, hence $\Orb_f(s)$ is finite, or $t$ is the fixed point $0$ of $g(x)$.

Let $q$ be an odd rational prime such that $q> r_i$ for all $i$. The latter argument shows that there is $\ell\mod q$ such that there are infinitely many pairs $(m_i,qi+\ell)$ where $f^{m_i}(s)=a^{qi+\ell}t$. This yields infinitely many $S$-integer points $(x,y)=(f^{m_i-1}(s),a^i)$ on the curve $a^{\ell}ty^q=f(x)$. 

Building on earlier work of Siegel, \cite{Siegel1,Siegel2}, Lang and LeVeque proved that if the number of $S$-integer points on a curve $C:y^{q}=f(x)$, $q\ge2$, $f(x)\in K[x]$, is infinite, then the genus of the curve $C$ must be zero, \cite{Lang,Veque1}. The reader may also consult \cite{Bu} for further references and literature. LeVeque also gave necessary and sufficient conditions for the genus of $C$ to be zero, \cite{Veque}. More precisely, setting $q_i=q/\gcd(q,r_i)$, and  assuming without loss of generality that $q_1\ge q_2\ge\cdots\ge q_e$, the curve $C$ has infinitely many $S$-integer points if and only if $(q_1,q_2,q_3,\cdots,q_e)=(2,2,1,\cdots,1)$ or $(s,1,1,\cdots,1)$, $s\ge1$. 

In view of the latter fact, since the tuple $(2,2,1,\cdots,1)$ is not realized due to the fact that $q$ is odd, either $f(x)$ has a root of multiplicity divisible by $q$, corresponding to the case  $e\geq 2$ and $q_2=1$, contradicting the assumption that $q>r_i$ for all $i$; or  $e=1$ implying that $f(x)$ is a power of a linear polynomial and concluding the proof.
\end{Proof}

\begin{Remark}\label{commoniterate}
The remark following Proposition 5.3 of \cite{Zieve1} shows that if $f$ and $g$ are non-monic linear polynomials such that $\Orb_f(s)\cap\Orb_g(t)$ is infinite, then $f$ and $g$ must have a common iterate.    
\end{Remark}

Proposition \ref{prop:linear} justifies the fact that we will only consider intersections of orbits of polynomials $f$ of arbitrary degrees with orbits of monic linear polynomials. In fact, we will mainly focus on the latter intersection when $f$ has integer coefficients. This may be justified by the following example.

\begin{Example}
Let $f(x)=x^d+a/b$, $d\ge 2$, $|b|>1$, $\gcd(a,b)=1$. One may easily see that $f^n(0)=a_n/b^{d^{n}}$ for some sequence $a_n\in\Z$, $n\ge1$ with $\gcd(a_n,b)=1$. In addition, $\Orb_{x+r/s}(0)=\{nr/s:n\ge 1\}\in\frac{1}{s}\Z.$ Therefore, $\Orb_f(0)\cap\Orb_{x+r/s}(0)$ contains only finitely many points for any choice of $r/s$. 
\end{Example}

\begin{Lemma}
\label{lem1}
Let $g(x)=x+ \frac{a}{b}$ with $a,b\in \Z$. Then $\Orb_{g}(t)=\displaystyle \frac{1}{b}\Orb_{x+a}(bt)$ for any $t\in \Q$.
\end{Lemma}
\begin{Proof}
   This follows immediately by observing that $g^n(t)= \dfrac{bt+na}{b}$ for any $n\ge1$.
\end{Proof}

In view of Lemma \ref{lem1}, it is sufficient to focus on orbits of linear polynomials of the form $g(x)=x+a$ with $a\in \Z$.

\begin{Lemma}  
Let $g_i(x)=x+m_i$, $m_i\in\Z$, $i=1,2$ and $ \gcd(m_1,m_2)=1$. Set $g(x)=x+m_1m_2$. Then $\Orb_{g}(t)= \Orb_{g_1}(t) \cap \Orb_{g_2}(t)$ for any $t\in \Q$.
\end{Lemma}
\begin{Proof}
Let $k$ be a rational number such that $k\in \Orb_g(t)$. Then we have $k=t+rm_1m_2$ for some $r\in \Z$. Hence $k\in \Orb_{g_1}(t) \cap \Orb_{g_2}(t) $. Now assume $k\in \Orb_{g_1}(t) \cap \Orb_{g_2}(t)$. Then $k=t+s_1m_1=t+s_2m_2$ for some $s_1,s_2\in \Z$. Since $\gcd(m_1,m_2)=1$, $m_1|s_2 $ and $m_2|s_1$. Thus, there exists $s\in\Z$ such that $k=t+ sm_1m_2$, i.e., $k\in \Orb_g(t) $. 
\end{Proof}

In fact, one has the following result for intersections of polynomial orbits with orbits of monic linear polynomials.

  \begin{Proposition} \label{main}  
     Let $f(x)$ be a polynomial in $\Z[x]$ and $g(x)=x+a$ where $a$ is an integer. Let $t$ be an integer such that $n$ is the minimum positive integer for which $f^{n}(t)\in \Orb_{g}^\pm(t)$. Then 
$f^{k}(t)\in \Orb_{g}^\pm(t) $ if and only if $n|k.$
\end{Proposition}
\begin{Proof}
We assume that $a=p_1^{\alpha_1}p_2^{\alpha_2}\dots p_e^{\alpha_e}$ where $p_1,\cdots,p_e$ are distinct primes.
As $f^n(t)=t+am$ for some $m\in \Z$, one obtains
$f^n(t)\equiv t \mod p_i^{\alpha_i}$ for all $i$.

If $n|k$, then $f^k(t)\equiv t \mod p_i^{\alpha_i}$. Since the primes $p_i$ are distinct, one has $f^k(t)\equiv t \mod  a$. 

Now one assumes that $f^{k}(t)\in \Orb_{g}^\pm(t)$. This gives
$f^k(t) \equiv t \mod  a$, hence
$f^k(t) \equiv t \mod p_i^{\alpha_i}$ for all $i=1,\cdots,e.$
We set $n_i$ to be the exact period of $t$ for the image $\widetilde{f}(x)$ of $f(x)$ in $\in\Z/(p_i^{\alpha_i}\Z)[x]$. Letting $l$ be the least common multiple of all $n_i$, $i=1,\dots,e$, one observes that $l|k$ as each $n_i |k$. In addition, one sees that $f^l(t)\equiv t\mod a$. This must yield that $l=n$, hence $n|k$.
\end{Proof}

\section{Primitive divisors and intersections with linear orbits}

We recall the following definition.
\begin{Definition}
    For an integer sequence $a_n$, $n\ge1$,
    a positive integer (rational prime) $u\ge 2$ is said to be a {\em primitive divisor (primitive prime divisor)} of $a_m$, if $u|a_m$ and $u\nmid a_s$ for all $1\leq s< m$.
\end{Definition}

For example, letting $\{a_n:n\geq 1\}$ be a sequence in which
    \[a_1=1,\quad a_2=2,\quad a_3=3,\quad a_4=6,\]
    one sees that $2$ is a primitive (prime) divisor for $a_2$,
    $3$ is a primitive (prime) divisor for $a_3$.
    However, $a_4$ does not have any primitive prime divisors, but $6$ is a primitive divisor of $a_4$.

Given a polynomial $f(x)\in\Z[x]$ and $t\in \Z$, we will investigate the set of primitive divisors of the sequence $t_0=0$ and $t_n=f^n(t)-t$, $n\ge 1$. 

\begin{Lemma}\label{lem:primitivedivisorssequence}
Let $t\in\Z$ and $f(x)\in\Z[x]$. Define the sequence $t_0=0$ and $t_n=f^n(t)-t$, $n\ge 1$. There exists a monic linear polynomial $g(x)=x+a\in\Z[x]$, $a\not\in\{0,\pm1\}$, such that $\Orb_{f}(t) \cap \Orb_{g}^\pm(t)\ne\{t\}$ if and only if $a$ is a primitive divisor of $t_m$ for some $m\ge1$. In this case, $\Orb_{f}(t) \cap \Orb_{g}^\pm(t)$ is infinite and $\Orb_{f^m}(t)\subset \Orb_g^\pm(t)$.
\end{Lemma}
\begin{Proof}
If $g(x)=x+a\in\Z[x]$ is such that $\Orb_{f}(t) \cap \Orb_{g}^\pm(t)\ne\{t\}$, then this is equivalent to the fact that there are integers $m,n\ge 1$ such that $f^m(t)=t+na$. Assuming that $m$ is the smallest such positive integer, one sees that $a$ is a primitive divisor of $f^m(t)-t$. The infinitude of the intersection and inclusion follows from Proposition \ref{main}. 
\end{Proof}

\begin{Lemma}
\label{cor:contain}
Let $f(x)\in\Z[x]$ and $t\in \Z$.  Let $h(x):=(f)^{\phi^{-1}}(x) \in\Z[x]$, where $\phi(x)=ax+t$. If $a$ is a divisor of $f(t)-t$, then the primitive prime divisors of $f^{i}(t)-t$ that do not divide $a$ are the primitive prime divisors of $h^{i}(0)$ that do not divide $a$ for all $i\ge 1$. 
\end{Lemma}
\begin{Proof}
Since $a\mid (f(t)-t)$, one easily sees that the map $h(x)=(f(t+ax)-t)/a\in\Z[x].$

One also sees that $h^i(0)=(f^i(t)-t)/a$. So, a prime $p\nmid a$ divides $h^i(0)$ if and only $p|f^{i}(t)-t$.

\end{Proof}

\begin{Lemma}

\label{onecover}
 Let $f(x)\in\Z[x]$. For all but finitely many integers $t$, there exists an integer $a\not\in\{0,\pm 1\}$, such that $\Orb_{f}(t)\subset\Orb_g^\pm(t)$ where $g(x)=x+a$. 
\end{Lemma}
\begin{Proof}
One sees that there are only finitely many $t$ such that $f(t)-t\in\{0,\pm 1\}$. For any other integer $t$, $f(t)-t$ has a primitive divisor $a$. Now the result follows by Lemma \ref{lem:primitivedivisorssequence}.  
\end{Proof}

In fact, Lemma \ref{onecover} can be made stronger by the following Proposition.

\begin{Theorem}\label{thm:onecover2} 
Let $f(x)\in\Z[x]$, and $m\geq 1$. For all but finitely many $t$, there exists an integer $a\not\in\{0,\pm1\},$ such that $\Orb_f(t)\cap \Orb_g^\pm(t)=\Orb_{f^m}(t)$ where $g(x)=x+a$.
\end{Theorem}

\begin{Proof}
Let $\phi_n(x)\in\Z[x]$ be the $n$-th dynatomic polynomial. Then there are only finitely many $t$ such that $\phi_m(t)\in\{0,\pm1\}$. For any other $t$, there exists a prime $p$ that divides $\phi_m(t)$. Let $r:=\nu_p(f^m(t)-t)$ and $g(x)=x+p^r$. Clearly, by the periodicity of $t$ under the reduction of $f^m(x)$ mod $p^r$, we get that $\Orb_{f^m}(t)\subset\Orb_g^\pm(t)$.

 Letting $n$ be such that $p$ is a primitive divisor of $f^n(t)-t$, then by Theorem 4.5 of \cite{dynamicsbook}, $m=n,\ m=nc,$ or $ m=ncp^e$, where $c$ is a constant and $e\geq 1$. For the first case where $m=n$, it is trivial that $f^k(t)\in\Orb_g^\pm(t)$ if and only if $m|k$.

 For the other cases, one can write $f^m(t)-t=\underset{d|m}{\prod}\phi_d(t)$ Since $p|\phi_m(t)$, then it is easy to see that for $s<m$, $\nu_p(f^s(t)-t)=\nu_p(\underset{d|s}{\prod}\phi_d(t))\leq \underset{d\leq s}{\sum}\nu_p(\phi_d(t))<r$ implying that $f^s(t)\not\in\Orb_g^\pm(t)$ and that $m$ is the least integer such that $f^m(t)\in\Orb_g^\pm(t)$. The rest follows from Proposition \ref{main}.
\end{Proof}

The latter results show that given a polynomial $f$ over $\Z$, there exists a linear polynomial $g$ that covers the orbit of $t$ under the iterates of $f$ for almost all integers $t$. Although a linear orbit cannot be covered by an orbit of a polynomial with higher degree, as will be shown in Lemma \ref{fun_facts_lemma}.

The next example shows that given a linear polynomial $g$ and a rational point $p$, there exists a quadratic polynomial $f$ such that the orbit of $p$ under $f$ is covered by the orbit of $p$ under $g$. In fact, one may establish the same result for a polynomial $f$ of any higher degree using an induction argument.

\begin{Example}
Consider the linear map $g(x)=x+k$ for some $k\in\Q^{\times}$. Let $a\in \Q^{\times}$ be such that $ak\in\Z$, and $p$ be a rational number.  There exists nonzero rational numbers $b,c$ such that the quadratic map $f(x)=ax^2+bx+c$ satisfies 
\[\Orb_f(p)\subset \Orb_g(p).\] More precisely,
\[\Orb_f(p)=\{p,p+k,p+2k,p+n_0k,p+n_1k,\cdots,p+n_i k,\cdots\}\]
where $n_0=3+2ak $, and $n_{i}= h^i(n_0)$, $i=1,2,\cdots$, where $h(x)= akx^2+(1-ak) x+1$.

For example, let $g(x)=x+7$, $f(x)= 2x^2-37x+163$, and $p=6$. Then one can have 
$$\Orb_f(p)=\{6,13,20,223,91370,\dots \} \subset \Orb_g(p)=\{6,13,20,26,\dots, 223,\dots, 91370, \dots \} .$$
\end{Example}

\section{Relative density of orbits intersections} \label{relativedens}
Given a polynomial $f$ in $\Z[x]$, an integer $s$, and a set $A\subseteq \Z$, we define the {\em relative density of $A$ in the orbit of $s$ under $f$} to be the limit
\[\delta_{f,s}(A):=\lim_{X\to\infty} \frac{| \{x\in A\cap\Orb_f(s):x\le X\}|}{| \{x\in \Orb_f(s):x\le X\}|},\] provided that this limit exists. 

Given a polynomial $g\in\Z[x]$ and an integer $t$, we will concern ourselves with $\delta_{f,s}(\Orb_g(t))$. The following lemma collects relative density results for some subsets of $\Z$ in the orbit of an integer under a polynomial $f$.

\begin{Lemma}\label{fun_facts_lemma} Let $f,g\in\Z[x]$ and $s,t\in\Z$ be such that $\Orb_f(s)$ is infinite. The following statements hold.
\begin{itemize}
\item[i)] $\delta_{f,s}(\Z)=1$.
\item[ii)] $\delta_{f,s}(\Orb_g(t))=0$ if $\deg f,\deg g>1$ where $f$ and $g$ have no common iterate; or $\deg f=\deg g=1$ where both $f$ and $g$ are non-monic and $f$ and $g$ have no common iterate.
\item[iii)] $\delta_{f,0}(\Orb_g(0))=1/m_2$, where $f(x)=x+m_1$, $g(x)=x+m_2$ and $\gcd(m_1,m_2)=1$. 
\item[iv)] $\delta_{f,s}(\Orb_g(t))=0$ if $\deg f=1$, $\deg g>1$ and $f$ is monic.
\end{itemize}
\end{Lemma}
\begin{Proof} 
For i), it is clear as $f\in\Z[x]$ and $t\in\Z$. 

For ii), if $\deg f,\deg g>1$ where $f$ and $g$ have no common iterate, Theorem 1 in  \cite{Zieve2} gives that $\Orb_{f}(s) \cap \Orb_{g}(t)$ is finite. Hence,  $\delta_{f,s}(\Orb_g(t))=0$. If $\deg f=\deg g=1$ where $f$ and $g$ have no common iterate and both $f$ and $g$ are non-monic, then the result follows from Remark \ref{commoniterate}.

For iii), the statement follows directly by observing that  $f^n(0)\in\Orb_g(0)$ if and only if $nm_1=km_2$ for some $k\in\Z$, and the coprimality of $m_1$ and $m_2$.

  For iv), we set $f(x)=x+b$, $b\in\Z$. We may assume without loss of generality that $t\ne0,\pm1$, otherwise, we may consider $g(t)$ instead. Let $m_1\geq 0$ be the least integer such that $g^{m_1}(t)\in \Orb_f(s)$, i.e., $g^{m_1}(t)=s+br_1$ for some $r_1 \in \Z$. Let $m_2$ be the least integer such that $m_2>m_1$ and  $g^{m_2}(t)\in \Orb_f(s)$, in particular, $g^{m_2}(t)=s+br_2$ for some $r_2\in \Z$. One can construct a sequence $m_1<m_2<\dots<m_i<\dots$ such that  $g^{m_i}(t)=s+br_i$ for some $r_i\in \Z$. We observe that 
$g^{m_i}(t)=g^{m_1}(t)+b(r_i-r_1)$, $i\ge 1$, hence $r_{i+1}-r_i=(g^{m_{i+1}}(t)-g^{m_i}(t))/b$. In other words, $r_{i+1}-r_i$ is a polynomial expression in $t$ of degree $d^{m_{i+1}}$, where $d=\deg g>1$. Since $m_{i}\geq i$, then $r_i\geq c\cdot|t|^{d^i}$ for some constant $c\in\Q$. Now, 

\begin{eqnarray*}\delta_{f,s}(\Orb_g(t))&=&\lim_{X\to\infty}\frac{|\{n\leq X:f^n(s)\in\Orb_g(t)\}|}{X}=\lim_{X\to\infty}\frac{|\{i:r_i\leq X\}|}{X}\\
&=&\lim_{X\to\infty}\frac{\max_{r_i\leq X} i}{X}=\lim_{i\to\infty}\frac{i}{r_i}\leq\lim_{i\to\infty}\frac{i}{c|t|^{d^i}}=0.
\end{eqnarray*}
\end{Proof}

\begin{Remark}
In Lemma \ref{fun_facts_lemma} iv), if $f(x)=ax+b$ where $a\ne \pm 1$, then Proposition \ref{prop:linear} implies that $\Orb_{f}^\pm(s) \cap \Orb_{g}(t)$ is finite, hence $\delta_{f,s}(\Orb_g(t))=0$. In addition, if $f(x)=-x+b$, then $\Orb_f^\pm(s)=\{s,b-s\}$, hence the orbit is finite in this case.
\end{Remark}

\begin{Corollary}\label{cor}
Let $f(x)$ be a polynomial in $\Z[x]$ of degree $d\ge 2$ and $g(x)=x+a$ where $a\in \Z$ is such that $a\ne0,\pm1$. Let $t\in \Z$. The following statements are equivalent:
\begin{itemize}
    \item[(i)] $t$ is a periodic point of $f(x) \mod a$ with exact period $n\ge1$.
    \item [(ii)] $a$ is primitive divisor of $f^n(t)-t$.
    \item[(iii)] $\delta_{f,t}(\Orb_{g}^\pm(t))=\frac{1}{n}$.
        \item[(iv)] $|\Orb_{f}(t) \cap \Orb_{g}^{\pm}(t)| = \infty$.
\end{itemize}
\end{Corollary}
\begin{Proof}
  This is Proposition \ref{main} and Lemma \ref{lem:primitivedivisorssequence}. 
\end{Proof}

In view of Corollary \ref{cor}, the infinitude of the intersection of a linear orbit with the orbit of an integer $t$ under a polynomial $f$ of arbitrary degree is equivalent to the existence of a primitive divisor for an element in the sequence  $\{f^{i}(t)-t\}_i$.

In fact, for the polynomial $f(x)=x^d+c\in\Z[x]$, it was proved in \cite{DoerksenHaensch+2012+465+472} that the sequence $f^i(0)$ has a primitive prime divisor for all but finitely many $i$'s. For $x^d+c\in\Q[x]$, it was shown in \cite{Holly_Krieger} that the sequence $f^i(0)$ has a primitive prime divisor for all $i$ except possibly for $23$ values. Moreover, it was shown in \cite{Rice} that for two classes of polynomials $f(x)\in \Z[x]$ and any integer $t$, the sequence $f^n(t)$, $n \geq 1$, has only finitely
many terms with no primitive prime divisor. 

If $t$ is a point whose orbit is infinite under $f(x)\in\Z[x]$, then for all but finitely many integers $n$, $f^n(t)$ has a primitive prime divisor under the $abc$-conjecture, see \cite{abc}. Moreover, if $t$ is a critical point of $f(x)$, then for all but finitely many integers $n$, $f^n(t)-t$ has a primitive prime divisor, see \cite{Ren}.

\begin{Remark}
In Corollary \ref{cor}, if $a=p_1^{\alpha_1}\cdots p_e^{\alpha_e}$, then $n$ is the least common multiple of all $n_i$, $i=1,\cdots,e$, where $n_i$ is the exact period of $t$ for the image $\widetilde{f}(x)$ of $f(x)$ in $\in\Z/(p_i^{\alpha_i}\Z)[x]$, see the proof of Proposition \ref{main}.
\end{Remark}

\section{Covering polynomial orbits using arithmetic progressions}
\label{sec5}
In this section, we discuss covers of polynomial orbits.

\begin{Definition}
Let $f(x)\in\Z[x]$ and $t$ be an integer. A finite system $$A=\{a_s+n_s\Z\}_{s=1}^k,\qquad a_s,n_s\in\Z,\quad n_s>0,\quad 1\le s\le k, $$ is said to be a {\em cover} of $\Orb_f(t)$ if $$ \Orb_f(t) \subset \cup_{s=1}^k \{a_s+n_s\Z\}.$$ If $\displaystyle \left\{a_s+n_s\Z\right\}_{\substack{s=1\\ s\ne u}}^k$ is not a cover of $\Orb_f(t)$, then $A$ is a cover of $\Orb_f(t)$ for which $a_u+n_u\Z$ is essential. A {\em minimal cover} of $\Orb_f(t)$ is a cover in which all the arithmetic sequences are essential. If $\Orb_f(t)\cap \{a_s+n_s\Z\}\cap \{a_{s'}+n_{s'}\Z\}=\emptyset$ for all $s\neq s'$, then $A$ is called a {\em disjoint cover} of $\Orb_f(t)$. A $t$-cover of $A$ is a cover of $\Orb_f(t)$ for which $t\in \cap_{s=1}^k \{a_s+n_s\Z\}$. 
\end{Definition}
For a positive integer $n$, one sees that $\{r + n\Z \}_{r=0}^{n-1}$ is a  disjoint cover of $\Orb_f(t)$ for any $f(x)\in\Z[x]$ and any integer $t$.

Given $f(x)\in\Z[x]$ and $t\in\Z$, we will mainly focus on $t$-covers of $\Orb_f(t)$ of the form $A=\{t+n_s\Z\}_{s=1}^k$. In other words, the latter $t$-covers can be expressed in the form 
$\{\Orb^\pm_{x+n_{s}}(t)\}_{s=1}^k$.

\begin{Theorem}\label{thm:cover}
 Let $f(x)$ be a polynomial in $\Z[x]$ of degree $d\ge 2$ and $t\in\Z$ be a wandering point for $f$. Let $g_i(x)$, $1\le i\le k$, be a finite family of monic linear polynomials in $\Z[x]$. The following statements are equivalent.
\begin{itemize}
    \item [i)] $\delta_{f,t} \left(  \bigcup_{i=1}^k \Orb_{g_i}^\pm(t) \right) = 1$.
    \item[ii)] $\delta_{f,t} \left(\Orb_{g_i}^\pm(t) \right) = 1$ for some $i$, $1\le i\le k$.
    \item[iii)]  $ \Orb_{f}(t) \subset  \Orb_{g_i}^\pm(t)$ for some $i$, $1\le i\le k$. 
\end{itemize}  
\end{Theorem}
\begin{Proof}
That the implication iii) yields i) is clear. 
 We assume that $\delta_{f,t} \left( \bigcup_i \Orb_{g_i}^\pm(t) \right) = 1$. We assume without loss of generality that $|\Orb_f(t) \cap \Orb_{g_{i}}^\pm(t)|= \infty $ for all $i$. By Corollary \ref{cor}, $\delta_{f,t} ( \Orb_{g_{i}}^\pm(t))=\frac{1}{n_i}$ for some positive integer $n_i$. Moreover, if $g_i(x)=x+a_i$, then $a_i$ is a primitive divisor of $f^{n_i}(t)-t$. We assume that $n_i > 1$ for all $i$. Proposition \ref{main} implies that
 $f^{m}(t) \in \Orb_{g_i}^\pm(t)$ if and only if $n_i|m$. Setting $n=\prod_i n_i$, it is obvious that $n_i\nmid (hn+1)$ for any integer $h$. In particular,
  $f^{hn+1}(t) \notin \Orb_{g_i}^\pm(t)$ for any $i$, $1\le i\le k$. This implies that
 $\delta_{f,t} \left( \bigcup_i \Orb_{g_i}^\pm(t) \right) \leq 1-\frac{1}{n};$ which contradicts our assumption, hence $n_i=1$ for some $i$. 
 
 Assuming ii), Corollary \ref{main} ii) implies that $a_i$ is a primitive divisor of $f(t)-t$. In view of Lemma \ref{lem:primitivedivisorssequence}, one sees that
 $\Orb_f(t) \subset \Orb_{g_i}^\pm(t)$.
 \end{Proof}
\begin{Corollary}
Let $f(x)$ be a polynomial in $\Z[x]$ of degree $d\ge 2$ and $t\in\Z$. If $A=\{t+n_s\Z\}_{s=1}^k$ is a minimal $t$-cover of $\Orb_f(t)$, then $k=1$. In particular, if $A=\{t+n_s\Z\}_{s=1}^k$, $k\ge 1$, is a system of arithmetic progressions such that $\delta_{f,t}\left(t+n_s\Z\right)<1$, for all $s=1,\cdots,k$, then $\delta_{f,t}\left(\bigcup_{i=1}^k(t+n_s\Z)\right)<1$.
\end{Corollary}

\section{Densities of intersections of polynomial orbits and arithmetic progressions}

In \S\ref{sec5}, we showed that a polynomial orbit of a rational point $t$ cannot be covered by finitely many arithmetic progressions, each of which is containing $t$, unless the orbit lies in one of these progressions. In this section, we study the density of the intersection of the orbit with finitely many arithmetic progressions, each of which is containing $t$. 
 We investigate the set of real numbers that appear as such densities. Examining the properties of these numbers, we bound the density of the intersection from above using a bound that depends only on the number of arithmetic progressions. We also show that this upper bound can be arbitrarily close to $1$ if the number of progressions increases.

\begin{Definition}
Let $f(x)\in\Z[x]$ and $t\in\Z$. Let $\delta$ be a real number such that $0\le \delta\le 1$. If there is a system of arithmetic progressions of the form  $A=\{t+n_s\Z\}_{s=1}^k$ such that $\delta_{f,t}\left(\bigcup_{s=1}^kt+n_s\Z\right)=\delta$, then $\delta$ is said to be {\em $(f,t,k)$-accessible}.   
\end{Definition}

For $f(x)\in\Z[x]$ and $t\in\Z$, we set $$\mathcal{PD}(f,t)=\{a: a \textrm{ is a primitive divisor of } f^n(t)-t \textrm{ for some }n\ge 1\}.$$ The set $\mathcal{PD}(f,t)$ contains the set of primitive prime divisors of $f^n(t)-t$. If $t$ is a wandering point for $f$, then it is clear that $\mathcal{PD}(f,t)$ is infinite since otherwise $t$ will be a preperiodic point under $f$. We also set 
$$S(f,t)=\{n\ge 1: \,f^n(t)-t \textrm{ has a primitive divisor}\}.$$
Again, if $t$ is a wandering point for $f$, then $S(f,t)$ is infinite.

\begin{Definition}\label{inc-ex-def}
Let $S\subseteq\Z^{+}=\{z\in\Z:z>0\}$. A nonnegative rational number $\delta<1$ is said to be an $(S,k)$-{\em inclusion-exclusion fraction} if there are $n_i\in S$, $i=1,\cdots,k$, with
$$\delta=\sum_{i=1}^k \frac{1}{n_i}-\sum_{1\le i_1<i_2\le k}\frac{1}{\lcm(n_{i_1},n_{i_2})}+\sum_{1\le i_1<i_2<i_3\le k}\frac{1}{\lcm(n_{i_1},n_{i_2},n_{i_3})}+\cdots+(-1)^{k+1}\frac{1}{\lcm(n_1,\cdots, n_k)}.$$
\end{Definition}

\begin{Theorem}\label{accessible}
Let $f(x)\in\Z[x]$ be of degree $d\ge 2$ and $t\in\Z$ be a wandering point for $f$. Let $k\ge1$ be an integer. 
An $(S(f,t),k)$-inclusion-exclusion fraction is $(f,t,k)$-accessible.
\end{Theorem}
\begin{Proof}
   Let $\delta$ be an $(S(f,t),k)$-inclusion-exclusion fraction where $n_i\in S(f,t)$, $1\leq i\leq k$, are as in Definition \ref{inc-ex-def}. Let $a_i$ be a primitive divisor of $f^{n_i}(t)-t$, see the definition of $S(f,t)$. Let $g_i=x+a_i$. By Corollary \ref{cor}, $\delta_{f,t}(\Orb_{g_i}^\pm(t))=\frac{1}{n_i}$. By setting $A=\{t+a_i\Z\}$, we can see that

\begin{eqnarray*}
    \delta_{f,t}\left(\bigcup_{i=1}^kt+a_i\Z\right)&=&\delta_{f,t}\left(\bigcup_{i=1}^k\Orb_{g_i}^\pm(t)\right)\\
    &=&\lim_{X\to\infty}\frac{|\{x\in \bigcup_{i=1}^{k}\left(\Orb_{g_i}^\pm(t)\cap \Orb_f(t)\right):x\leq X\}|}{|x\in  \Orb_f(t):x\leq X\}|}\\
   &=&\sum_{j=1}^{k}(-1)^{j+1}\sum_{1\leq i_1<i_2<\dots<i_j\leq k}\lim_{X\to\infty}\frac{|\{x\in \bigcup_{r=1}^{j}\left(\Orb_{g_{i_r}}^\pm(t)\cap \Orb_f(t)\right):x\leq X\}|}{|x\in  \Orb_f(t):x\leq X\}|}\\
    &=&\sum_{j=1}^{k}(-1)^{j+1}\sum_{1\leq i_1<i_2<\dots<i_j\leq k}\lim_{M\to\infty}\frac{|\{f^m(t)\in \bigcup_{r=1}^{j}\Orb_{g_{i_r}}^\pm(t):m\leq M\}|}{M}\\
    &=&\sum_{j=1}^{k}(-1)^{j+1}\sum_{1\leq i_1<i_2<\dots<i_j\leq k}\frac{1}{\lcm(n_{i_1},\dots,n_{i_j})}\\
\end{eqnarray*}

which concludes the result. The third equality is by the inclusion-exclusion principle. The last equality follows by Corollary \ref{cor} and Proposition \ref{main}, since $f^m(t)\in \Orb_{g_{i_r}}^\pm(t)$ if and only if $n_{i_r}|m$ which implies that $f^m(t)\in \bigcap_{r=1}^{j}\Orb_{g_{i_r}}^\pm(t)$ if and only if $\lcm(n_{i_1},\dots,n_{i_j})|m$.

\end{Proof}

\begin{Corollary}\label{1/n accessible}
    Let $f(x)\in\Z[x]$ be of degree $d\ge 2$ and $t\in\Z$ be a wandering point for $f$. If $n\in S(f,t)$, then $1/n$ is $(f,t,1)$-accessible.  
\end{Corollary}
\begin{Proof}
    This follows from Corollary \ref{cor}.
\end{Proof}

\begin{Proposition}\label{2/n}
Let $f(x)\in\Z[x]$ be of degree $d\ge 2$ and $t\in\Z$ be a wandering point for $f$. 
\begin{itemize}
\item[i)] If $n,\, n-1\in S(f,t)$, then $2/n$ is $(f,t,2)$-accessible.
\item[ii)] Let $n$ be an odd integer such that $n,\, n-1,\,n-2\in S(f,t)$, then $3/n$ is $(f,t,3)$-accessible.
\item[iii)]  Let $m,n\in\Z$ such that $(m-1)|(n-1)$ and $n, (n-1)/(m-1) \in S(f,t)$, then $m/n$ is $(f,t,2)$-accessible.
   \end{itemize}
\end{Proposition}

\begin{Proof} 
For i), if $n,\, n-1\in S(f,t)$ then  $1/n$ and $1/(n-1)$ are $(f,t,1)$-accessible from Corollary \ref{1/n accessible}. Since $\gcd(n,n-1)=1$, one can have that
$$\frac{1}{n}+\frac{1}{n-1}- \frac{1}{n(n-1)}=\frac{2}{n}$$
is $(f,t,2)$-accessible by using Theorem \ref{accessible}.

For ii), Corollary \ref{1/n accessible} and assumption give that $1/n, 1/(n-1)$ and  $1/(n-2)$ are $(f,t,1)$-accessible. Since $n$ is odd integer, we have $\gcd(n,n-1)=\gcd(n,n-2)=\gcd(n-1,n-2)=1$. Then by using Theorem \ref{1/n accessible}, one sees that
$$\frac{1}{n}+\frac{1}{n-1}+\frac{1}{n-2}-\frac{1}{n(n-1)}-\frac{1}{n(n-2)}-\frac{1}{(n-1)(n-2)}+\frac{1}{n(n-1)(n-2)}= \frac{3}{n}$$
is $(f,t,3)$-accessible.

The proof of iii) is similar. 

\end{Proof}

Given $f(x)\in\Z[x]$ and a wandering point $t\in\Z$ for $f$, Corollary \ref{1/n accessible} and Proposition \ref{2/n} display some $(f,t,k)$-accessible rationals, $k\ge 1$. These rationals are arising as inclusion-exclusions fractions. In what follows, fixing an integer $k\ge 1$, we introduce rationals in the interval $(0,1)$ that are not $(f,t,k)$-accessible.

   \begin{Theorem}\label{notexist}
       Let $f(x)\in\Z[x]$ be of degree $d\ge 2$ and $t\in\Z$ be a wandering point for $f$. Let $r\in (0,1)$ and $p_i$ be the sequence of ordered rational primes. If $k$ is a positive integer such that  $$\delta_k= 1- \prod_{i=1}^{k}\left(1-\frac{1}{p_i} \right) < r, $$
      then $r$ is not $(f,t,k)$-accessible. In particular, there does not exist $k$ linear polynomials $g_1(x),\dots,g_k(x)$ such that
     $$\delta_{f,t} \left(  \bigcup_{i=1}^{k} \Orb_{g_i}^\pm(t) \right)= r.$$
   \end{Theorem} 

\begin{Proof}
 In what follows we will show that 

$$S_k:=\underset{{(g_1,\cdots,g_k)\in \{x+a:a\in \Z\}^k}}{\sup}{\delta_{f,t}\left(  \bigcup_{i=1}^{k} \Orb_{g_i}^\pm(t) \right)}\le \delta_k$$ which immediately proves the statement.

We start by recalling that if $g(x)=x+a$, $a\in \Z$, is such that $|\Orb_f(t)\cap\Orb_g^\pm(t)|=\infty$, then $\Orb_f(t)\cap\Orb_g^\pm(t)=\Orb_{f^n}(t)$ for some $n\ge 1$, see Corollary \ref{cor}. 
    Now, given $g_1(x)=x+a_1$ and $g_2(x)=x+a_2$ with $\Orb_f(t)\cap\Orb_{g_i}^\pm(t)=\Orb_{f^{n_i}}(t)$ such that $\gcd(n_1,n_2)=n_{12}$, one notices that $\Orb_{f^{n_1}}(t)\cup\Orb_{f^{n_2}}(t)\subset\Orb_{f^{n_{12}}}(t)$. In particular, one only needs to examine the tuples $(g_1,\cdots,g_k)$ of monic linear polynomials for which $\Orb_f(t)\cap\Orb_{g_i}^\pm(t)=\Orb_{f^{n_i}}(t)$ is such that $n_1,\cdots,n_k$ are pairwise relatively prime. The following identities now hold 
    \begin{eqnarray*}
        \delta_{\{n_i\}_{i=1}^{k}}&:=&\sum_{i=1}^k \frac{1}{n_i}-\sum_{1\le i_1<i_2\le k}\frac{1}{\lcm(n_{i_1},n_{i_2})}+\sum_{1\le i_1<i_2<i_3\le k}\frac{1}{\lcm(n_{i_1},n_{i_2},n_{i_3})}+\cdots+(-1)^{k+1}\frac{1}{\lcm(n_1,\cdots, n_k)}\\
        &=&\sum_{i=1}^k \frac{1}{n_i}-\sum_{1\le i_1<i_2\le k}\frac{1}{n_{i_1}n_{i_2}}+\sum_{1\le i_1<i_2<i_3\le k}\frac{1}{n_{i_1}n_{i_2}n_{i_3}}+\cdots+(-1)^{k+1}\frac{1}{n_1\cdots n_k}\\
        &=&\sum_{i=2}^k \frac{1}{n_i}-\sum_{2\le i_1<i_2\le k}\frac{1}{n_{i_1}n_{i_2}}+\sum_{2\le i_1<i_2<i_3\le k}\frac{1}{n_{i_1}n_{i_2}n_{i_3}}+\cdots+(-1)^{k}\frac{1}{n_2\cdots n_k}\\
        &+&\frac{1}{n_1}\bigg(1-\sum_{i=2}^k \frac{1}{n_i}+\sum_{2\le i_1<i_2\le k}\frac{1}{n_{i_1}n_{i_2}}-\sum_{2\le i_1<i_2<i_3\le k}\frac{1}{n_{i_1}n_{i_2}n_{i_3}}-\cdots-(-1)^{k}\frac{1}{n_2\cdots n_k}\bigg)\\
        &=&\delta_{\{n_i\}_{i=2}^k}+\frac{1}{n_1}\left(1-\delta_{\{n_i\}_{i=2}^k}\right).
        \end{eqnarray*}
        More generally, after relabelling one sees that $\delta_{\{n_i\}_{i=1}^{k}}=\delta_{\{n_i\}_{i\ne j}}+\frac{1}{n_j}\left(1-\delta_{\{n_i\}_{i\ne j}}\right)$. Using an induction argument, it follows that 
        $\delta_{\{n_i\}_{i=1}^{k}}\le \delta_{\{q_i\}_{i=1}^k}$, where $q_i$ is a prime divisor of $n_i$. Now it follows that $S_k\le \delta_{\{p_i\}_{i=1}^k}$,
        where direct computations show that $\delta_{\{p_i\}_{i=1}^k}=1- \prod_{i=1}^{k}\left(1-\frac{1}{p_i} \right)$.
\end{Proof}

Given $f(x)\in\Z[x]$ together with a wandering point $t\in\Z$, Corollary \ref{1/n accessible} and Proposition \ref{2/n} assert that $i/n$, $i=1,2,3$, is $(f,t,i)$-accessible under certain conditions on $n$. Theorem \ref{notexist} shows that this phenomenon does not hold in general. 

\begin{Example}

Let $f(x)$ be a polynomial in $\Z[x]$ of degree $d\ge 2$ and $t\in\Z$. There are no linear polynomials  $g_1(x),\dots,g_4(x) \in \Z[x]$  such that
    $$\delta_{f,t}(\Orb_{g_1}^\pm(t) \cup \Orb_{g_2}^\pm(t) \cup \Orb_{g_3}^\pm(t)  \cup \Orb_{g_4}^\pm(t))=\frac{4}{5}.$$
    This holds because $\delta_4=27/35 <4/5$, see Theorem \ref{notexist}. 
    In particular, if $\delta_{f,t}\left(\cup_{i=1}^k\Orb_{g_i}^\pm(t)\right)=4/5$, then $k\ge 5$.   
\end{Example}
The following table summarizes the values in $r\in(0,1)$ that are not $(f,t,k)$-accessible for any $f\in\Z[x]$ with a wandering point $t\in \Z$, for $k\le 10$.
\begin{center}
{\footnotesize\begin{tabular}{| c | c | c|c|c|c|c|c|c|c|c|}
\hline
$k$ & $1$ &$2$&$3$&$4$&$5$&$6$&$7$& $8$&$9$ & $10$\\ 
\hline
$r>$ & $1/2$ & $2/3$ & $11/15$ & $27/35$ & $61/77$ & $809/1001$& $13945/17017$ & $268027/323323$ & $565447/676039$ & $2358365/2800733$\\
\hline
\end{tabular}}
\end{center}

In \cite[Theorem 1.2]{Ren}, it was proved that the set of integers $n$ such that $f^n(t)-t$ has a primitive prime divisor contains all but finitely many integers under the condition that $t$ is a critical point for $f$. In fact, the latter fact was proved for a polynomial with rational coefficients of degree at least $2$. For our purpose, we want to avoid any other restrictions on $t$ but being a wandering point for $f$. Therefore, we prove the following result.

\begin{Proposition}\label{Prop:P}
Let $f(x)$ be a polynomial in $\Z[x]$ of degree $d\ge 2$ and $t\in\Z$ be a wandering point for $f$. The set $S(f,t)$ contains all rational primes.  
\end{Proposition}
\begin{Proof}
Let $p$ be a rational prime. Since $t$ is a wandering point, it follows that $f^p(t)\neq f(t)$, in particular, $f^p(t)-t\neq f(t)-t$. We now consider the prime factorization $f(t)-t=\underset{i=1}{\overset{s}{\prod}}q_i^{\alpha_i}$ and $f^p(t)-t=\underset{i=1}{\overset{s}{\prod}}q_i^{\beta_i}$, where $q_i\ne q_j$ for $i\ne j$, and $\alpha_i,\beta_i\geq 0$ are such that $\alpha_i+\beta_i\geq 1$. 

We note that since $f(t)\equiv t \mod q_i^{\alpha_i}$, we have $f^p(t)\equiv t \mod q_i^{\alpha_i}$, $i=1,\cdots,s$. In particular, we see that $\alpha_i\leq \beta_i$, $i=1,\cdots,s$. Since $f^p(t)-t\neq f(t)-t$, it follows that there exists at least one $1\leq j\leq s$ such that $\alpha_j<\beta_j$. Noting that $f^p(t)\equiv t \mod q_j^{\beta_j}$, we obtain that $t$ is a periodic point of $f(x)$ mod $q_j^{\beta_j}$, where the exact period of $t$ divides $p$. Since $\beta_j>\alpha_j$, this yields that $p$ is the exact period for $t$. In view of Corollary \ref{cor}, we get that $q_j^{\beta_j}$ is a primitive divisor of $f^p(t)-t$, hence $p$ lies in $S(f,t)$.
\end{Proof}

Given a positive integer $k$, Theorem \ref{notexist} asserts that there is a rational number $a$, $0<a<1$, such that for any $r\in (a,1)$, $r$ is not $(f,t,k)$-accessible for any $f\in\Z[x]$ and $t\in\Z$.
The following Theorem shows that the interval $(a,1)$ can be made arbitrarily small for large enough values of $k$. We recall that the latter interval can never be empty due to Theorem \ref{thm:cover}.

\begin{Theorem}
    Let $f(x)$ be a polynomial in $\Z[x]$ of degree at least $2$, and $t\in\Z$ be a wandering point for $f$. Given an $\epsilon >0$, there exists finitely many linear polynomials $g_{i}(x)$, $i=1,\ldots,k(\epsilon)>1$, such that 
    $$1-\epsilon\le \delta_{f,t} \left( \bigcup_{i=1}^{k(\epsilon)} \Orb_{g_i}^\pm(t) \right)<1. $$
\end{Theorem}

\begin{Proof}
According to Proposition \ref{Prop:P}, $S(f,t)$ contains all rational primes $p_i$. 
Let $p_i$ be the exact period of $t$ under $f(x) \mod m_i$ for some positive integer $m_i$, or equivalently, $m_i$ is a primitive divisor of $f^{p_i}(t)-t$, see Corollary \ref{cor}. We observe that 

    $$1- \prod_{r=1}^M\left(1-\frac{1}{p_i} \right) =\sum_{i} \frac{1}{p_i}- \sum_{i,j} \frac{1}{p_{i}p_{j}} + \sum_{i,j,k} \frac{1}{p_{i}p_{j}p_{k}}- \dots+(-1)^{M+1}\frac{1}{p_1\cdots p_M} =\delta_{f,t} \left( \bigcup_{r=1}^M \Orb_{x+m_r}^\pm(t) \right).$$
Thus, we obtain that
$\displaystyle\lim_{M \xrightarrow[]{}{\infty}}\delta_{f,t} \left( \bigcup_{r=1}^M \Orb_{g_r}^\pm(t) \right)=1  $, where $g_r(x)=x+m_r$.
Now the statement follows from the definition of the limit.

\end{Proof}

The theorem above shows that if $a\in (0,1)$ is a real number, then there exists an integer $k\ge 1$, and $b\in \Q\cap (a,1)$ such that $b$ is $(f,t,k)$-accessible.


\begin{thebibliography}{MM}
	\frenchspacing
	\renewcommand{\baselinestretch}{1}





\bibitem{Bell}
J. P. Bell, D. Ghioca and T. J. Tucker, 
The dynamical Mordell-Lang conjecture, 
{\em Mathematical Surveys and Monographs}, 210, American
Mathematical Society, Providence, RI, 2016. xiii+280 pp.


\bibitem{Bu}
Y. Bugeaud,
Bounds for the solutions of superelliptic equations,
{\em Compositio Mathematica}, 107 (1997): 187--19.


\bibitem{Burcu}
B. Barsak\c{c}i and M. Sadek,
Simultaneous rational periodic points of degree-2 rational maps,
{\em Journal of Number Theory}, 243 (2023), 715--728.

\bibitem{DoerksenHaensch+2012+465+472}
K. Doerksen and A. Haensch,
Primitive prime divisors in zero orbits of polynomials,
 {\em Integers},
12 (2012),
465--472. 

\bibitem{Erdos}
P. Erd\H{o}s,
On integers of the form $2^k + p$ and some related problems, 
{\em Summa Brasil. Math.}, 2 (1950), 113--123.



\bibitem{Zieve1}
D. Ghioca, T. Tucker and M. Zieve,
Intersections of polynomial orbits, and a dynamical Mordell-Lang conjecture, 
{\em Inventiones Math.}, 
171 (2008), 463--483.

\bibitem{Zieve2}
D. Ghioca, T. Tucker and M. Zieve,
Linear relations between polynomial orbits, 
{\em Duke Math. J.}, 161 (2012), 1379--1410.

\bibitem{abc}
C. Gratton, K. Nguyen and T. J. Tucker,
ABC implies primitive prime divisors in arithmetic dynamics,
{\em Bulletin of the London Mathematical Society}, 45 (2013), 1194--1208.

\bibitem{Rafe_jones}
 S. Hamblen, R. Jones and K. Madhu,
 The density of primes in orbits of $z^d+c$,
  {\em International Mathematics Research Notices}, 
7 (2015), 1924--1958.

\bibitem{Hindes}
W. Hindes,
 Finite orbit points for sets of quadratic polynomials, 
{\em Int. J. Number Theory}, 15.8 (2019), 1693--1719.

\bibitem{Ho}
B. Hough,
Solution of the minimum modulus problem for covering systems, 
{\em Annals of Mathematics}, 181 (2015), 361--382.

\bibitem{Holly_Krieger}
H. Krieger,
Primitive prime divisors in the critical orbit of $z^d+c$,
 {\em International Mathematics Research Notices}, 
 23 (2013), 5498--5525.

\bibitem{Lang}
S. Lang, 
Integral points on curves, 
{\em Publ. Math. I.H.E.S.}, 6 (1960), 27--43.

\bibitem{Veque1} 
W. J. LeVeque, 
Rational points on curves of genus greater than 1,
{\em J. Reine Angew. Math.},
206 (1961), 45--52. 


\bibitem{Veque}
W. J. LeVeque, 
On the equation $y^m = f(x)$, 
{\em Acta Arith.}, 9 (1964), 209--219.


\bibitem{Ren}
R. Ren,
Primitive prime divisors in the critical orbits of one-parameter families
of rational polynomials,
{\em Math. Proc. Camb. Phil. Soc.}, 171 (2021), 569--584.


\bibitem{Rice}
B. Rice, 
Primitive Prime Divisors in Polynomial Arithmetic Dynamics,
{\em INTEGERS: Electronic Journal of Combinatorial Number Theory}, 7 (2007), Article 26.

\bibitem{SadekWafik}
M. Sadek and M. Wafik,
Construction of polynomials with prescribed divisibility conditions on the critical orbit, preprint.

\bibitem{SadekYesin}
M. Sadek and T. Yesin,
A Dynamical Analogue of a question of Fermat, accepted for publication in {\em Mathematica Scandinavica}.




\bibitem{Siegel1}
C. L. Siegel, 
The integer solutions of the equation $y^2 = ax^n + bx^{n-1} +\cdots+ k$,
{\em J. London Math. Soc.}, 1 (1926), 66--68; also: Gesammelte Abhandlungen, Band 1, 207--208.

\bibitem{Siegel2}
C. L. Siegel, 
\"{U}ber einige Anwendungen Diophantischer Approximationen, Abh. Preuss. Akad. Wiss. Phys.-Math. Kl., 1929, Nr. 1; also: Gesammelte Abhandlungen, Band 1, 209--266.


\bibitem{dynamicsbook}
 J. H. Silverman,
  The arithmetic of dynamical systems,
  Springer, 2007.

\end{thebibliography}
\end{document}